\documentclass[11 pt, psamsfonts]{amsart}
\usepackage{amssymb}

\def\ignore#1{\relax}

\def\g{\mathfrak g}

\def\o{\mathfrak o}

\def\R{{\mathbb R}}
\def\Z{{\mathbb Z}}

\def\C{{\mathbb C}}
\def\F{{\mathcal F}}
\def\la{\lambda}

\def\Ca{\mathcal C}

\def\E{\mathcal E}

\def\U{{\bf U}}

\def\ignore#1{\relax}

\def\om{\omega}

\def\1{{\bf 1}}

\def\ep{\epsilon}

\def\End{{\rm End}}

\def\ve{\varepsilon}

\calclayout

\renewenvironment{proof}[1][\proofname]{\par
  \normalfont
  \topsep6\p@\@plus6\p@ \trivlist
  \item[\hskip\labelsep\bfseries
    #1\@addpunct{.}]\ignorespaces
}{%
  \qed\endtrivlist
}

\def\th@plain{%
  \let\thmhead\thmhead@plain \let\swappedhead\swappedhead@plain
  \thm@preskip.5\baselineskip\@plus.2\baselineskip
                                    \@minus.2\baselineskip
  \thm@postskip\thm@preskip
  \itshape
\renewcommand{\labelenumi}{{(\alph{enumi})\quad}}
                        \renewcommand{\labelenumii}{{(\roman{enumii})\ }}
}
\def\th@definition{%
  \let\thmhead\thmhead@plain \let\swappedhead\swappedhead@plain
  \thm@preskip.5\baselineskip\@plus.2\baselineskip
                                    \@minus.2\baselineskip
  \thm@postskip\thm@preskip
  \upshape
}
\def\th@remark{%
  \thm@headfont{\itshape}
  \let\thmhead\thmhead@plain \let\swappedhead\swappedhead@plain
  \thm@preskip.5\baselineskip\@plus.2\baselineskip
                                    \@minus.2\baselineskip
  \thm@postskip\thm@preskip
  \upshape
}

{\theoremstyle{plain}
\newtheorem{theorem}{Theorem}[section]
}

{\theoremstyle{plain}
\newtheorem{proposition}[theorem]{Proposition}
}

{\theoremstyle{plain}
\newtheorem{corollary}[theorem]{Corollary}
}

{\theoremstyle{plain}
\newtheorem{lemma}[theorem]{Lemma}
}

{\theoremstyle{plain}

}

{\theoremstyle{definition}
\newtheorem{definition}[theorem]{Definition}
}

{\theoremstyle{definition}

}

{\theoremstyle{remark}
\newtheorem{remark}[theorem]{Remark}
}

{\theoremstyle{remark}

}

\numberwithin{equation}{section}

\renewcommand{\labelenumi}{{ \theenumi.}}
\renewcommand{\labelenumii}{{(\alph{enumii})}}

\def\la{\lambda}
\def\al{\alpha}

\def\so{{\mathfrak s}{\mathfrak o}}
\def\sp{{\mathfrak s}{\mathfrak p}}
\def\sl{{\mathfrak s}{\mathfrak l}}

\def\omi{\om_i}
\def\omj{\om_j}

\def\choose #1 #2{\begin{pmatrix}#1\\#2\end{pmatrix}}

\def\x{{\bf x}}

\ignore{
\input BoxedEPS
\SetRokickiEPSFSpecial
\HideDisplacementBoxes
\SetEPSFDirectory{./graphics/} }

\ignore{
\input BoxedEPS
\SetTexturesEPSFSpecial
\HideDisplacementBoxes
\SetEPSFDirectory{:graphics:}
}

\def\E{\mathcal E}
\def\Ca{\mathcal C}

\def\x{{\bf x}}

\def\U{{\bf U}}

\def\ignore#1{\relax}

\def\om{\omega}

\def\1{{\bf 1}}

\def\ep{\epsilon}

\def\End{{\rm End}}

\def\ve{\varepsilon}

\def\p{\psi}
\def\pd{\psi^\dagger}

\def\tve{\tilde\ve}
\def\bX{\bar X}

\def\Om{\Omega}

\begin{document}

\title[Dualities for Spin Representations ]
{Dualities for Spin Representations }

\author{Hans Wenzl}

\address{Department of Mathematics\\ University of California\\ San Diego,
California}

\email{hwenzl@ucsd.edu}

\begin{abstract}
Let $S$ be the spinor representation of $U_q\so_N$, for $N$ odd
and $q^2$ not a rooot of unity.
We show that the commutant of its action on $S^{\otimes n}$ is given
by a representation of the nonstandard quantum group $U'_{-q^2}\so_n$.
For $N$ even, an analogous statement also holds for $S=S_+\oplus S_-$
the direct sum of the irreducible spinor representations of $U'_q\so_N$,
with the commutant given by  $U'_{-q}\o_n$, 
a $\Z/2$-extension of $U'_{-q}\so_n$. Similar statements also hold
for fusion tensor categories with $q$ a root of unity.
\end{abstract}
\maketitle

The decomposition of tensor powers of the vector representations of
classical Lie groups was successfully studied in papers by Schur, Weyl and Brauer.
These classical  duality results were extended to Drinfeld-Jimbo quantum groups
(\cite{Ji}, \cite{BBMW}, \cite{WBMW})
which have had applications in a number of fields such as
tensor categories, low-dimensional topology and von Neumann algebras. But it seems that
an intrinsic description of the commutant of the action of spin groups 
on tensor powers of spinor representations has been found
only fairly recently in \cite{Wspin}. The current paper deals with
 the missing cases in that paper. Moreover,
we give  simpler and more general
proofs also for the cases already covered.

The inspiration for the approach in the paper \cite{Wspin} came from a
paper by Hasegawa \cite{Has}, see Section \ref{Hasegawa}
for more details. He showed (as a special case of a
far more general construction)
that the obviously commuting actions of $O(N)$ and $SO(n)$ on
$\R^N\otimes \R^n$ can be extended to commuting actions 
of the corresponding spin groups on the Clifford algebra
$Cl(Nn)\cong Cl(N)^{\otimes n}$. For $N$ even, this vector space
isomorphism can actually be made into an algebra homomorphism, with $Cl(N)\cong
\End(S)$ for the spinor representation $S$ of $O(N)$.
Hence the commutant of the $Pin(N)$ action on $S^{\otimes n}$
is given by a representation of $Spin(n)$. 
This result could be extended to prove a duality on $S^{\otimes n}$ between
a semidirect product $\U$  of the quantum group $U_q\so_N$ with $\Z/2$
and a non-standard $q$-deformation $U'_q\so_n$ of the universal
enveloping algebra $U\so_n$. The latter has been studied before
in particular by Klimyk and his coauthors, see e.g. \cite{GK}, and it has
also appeared in work of  Noumi and Sugitani \cite{NS} and Letzter \cite{Le}.

 But for $N$ odd, it was not possible to prove such a simple fact,
due to the fact that the Clifford algebra is not simple in that case.
 Our main new result in this paper
is a direct description of the commutant of the action of
$U_q\so_N$ on $S^{\otimes n}$ via a representation of
the nonstandard orthogonal quantum group $U'_{-q^2}\so_n$.
This can be quite easily seen for $q=1$ using the element
$C=\frac{1}{2}\sum_i e_i\otimes e_i\subset Cl(N)^{\otimes 2}$,
where $\{ e_i\}$ is an orthonormal basis for $\R^N$.
We then extend this result to the quantum group case
by using Hayashi's $q$-Clifford algebra.

Here is the paper in more detail: We first review some basics about
Clifford algebras and spin representations $S$ and produce
a canonical element $C\in\End(S^{\otimes 2})$ which commutes
with the spin group action. We review basic facts about the algebras
$U'_q\so_n$ and their representations in the second section.
In particular, we show that we obtain a representation of
$U'_{-1}\so_n$ on $S^{\otimes n}$ which commutes with the 
pin group action.
In the third section we review Hayashi's spin representations
of $U_q\so_N$
into his $q$-Clifford algebra and we construct the $q$-analogs of
the commuting elements $C$. It is shown in the following section
that these can be used to define a representation of $U'_{-q^2}\so_n$
on $S^{\otimes n}$. We then obtain first and second fundamental 
theorems for $\End_\U(S^{\otimes n})$, where $\U$
can be  $U_q\so_N$ or, for $N$ even, it can also be $U_q\so_N\rtimes \Z/2$.
We conclude with some remarks about related work and applications.

$Acknowledgements:$ Part of the work on this paper was done
while the author enjoyed the hospitality and support of the Max Planck Institute
of Mathematics in Bonn. He would like to thank Catharina Stroppel and
Daniel Tubbenhauer for stimulating discussions. He would also like to thank Doron Gepner 
for his interest and encouragement.

\section{Clifford Algebras and spinor representations}

\subsection{Basic Definitions}
Let $V$ be a finite dimensional real inner product space. Then its
Clifford algebra is the associative complex algebra generated by the elements of 
$V$ subject to the relation
$$ v\cdot w+w\cdot v=2(v,w)1.$$
Let $\{ e_1,\ ...\ e_N\}$ be an orthonormal basis of the 
inner product space $V$. Then the Clifford algebra
$Cl =Cl(N)$ corresponding to $V$ can also be defined via generators,
denoted by $e_i$ as well, and relations
$$e_ie_j+e_je_i=2\delta_{ij}1,\quad {\rm for}\ 1\leq i,j\leq N.$$
If $N=2k$ is even, it will be convenient to use a second presentation
in terms of generators $\p_j$ and $\pd_j$, $1\leq j\leq k$ defined by
\begin{equation}\label{psidef}
\p_j=\frac{1}{2}(e_{2j-1}+ie_{2j}),\hskip 3em
\pd_j=\frac{1}{2}(e_{2j-1}-ie_{2j}).
\end{equation}
Then a generator $\p_j$ anticommutes with all other generators except $\pd_j$
where we have
$$\p_j\pd_j+\pd_j\p_j=1,\quad 1\leq j\leq k.$$

\subsection{Spinor representations}\label{spinreps} Let us assume $N=2k$ to be even.
It is well-known that $Cl(N)$ has dimension $2^N$ and it is isomorphic
to $M_{2^{N/2}}$, where $M_d$ denotes the $d\times d$ matrices.
Let $S$ be a simple $Cl(N)$-module, of dimension $2^{N/2}$.
The action of an element $g$ in the orthogonal group $O(N)$
on $V$ induces an automorphism $\al_g$ on $Cl(N)$, for each $g\in O(N)$.
As any automorphism of the $d\times d$ matrices is inner,
we obtain a projective representation $g\mapsto u_g$
of $O(N)$ on the $Cl(N)$-module $S$.
This projective representation can be made into an honest
representation of the universal covering groups $Pin(N)$ of $O(N)$. 
By restriction, the module $S$ becomes a $Pin(N-1)$-module,
which we will denote by $\tilde S$. It decomposes into the direct sum
of two simple projective $O(N-1)$-modules $\tilde S_+\oplus \tilde S_-$.
These two modules are isomorphic as projective $SO(N-1)$-modules
and simple.
We need the following relations, which are easy to prove:

\begin{lemma}\label{fnlemma} Let $N=2k$ be even.
Define $f_r=(-i)^{r(r-1)/2}e_1e_2\ ...\ e_r\in Cl(N)$, for $1\leq r\leq N$. Then we have

(a) $f_re_i=(-1)^re_if_r$ for $i>r$ and $f_re_i=(-1)^{r-1}e_if_r$ for $i\leq r$,

(b) $f_r^2=1$ and $f_rf_s=(-1)^{r(s-r)}f_sf_r$ for $r<s$.

(c) $f_2=e_1e_2=(\p_1\pd_1-\pd_1\p_1)$ and $f_N=(-1)^{k-1}\prod_{j=1}^k (\p_j\pd_j-\pd_j\p_j)$.

(d) $\al_g(f_r)=det(g)f_r$ and $u_gf_r=det(g)f_ru_g$ for $g\in O(r)$, $1\leq r\leq N$.
\end{lemma}

$Proof.$ Parts (a) - (c) are straightforward. Part (d) is checked by an
explicit calculation for $N=2$. The same calculation also works for
$g\in O(N)$ which is the identity matrix except for a $2\times 2$ diagonal
block. As such matrices generate $O(N)$, the claim follows.

\subsection{Explicit description}\label{explicit}
We will  need a more explicit description of the spin module $S$. Observe that
$Cl= Cl_+Cl_-$, where $Cl_+$ and $Cl_-$ are the subalgebras generated by
$1$ and the elements $\pd_i$ for $Cl_+$, and by 1 and the elements $\p_i$
for $Cl_-$. We define $S=Cl/I$ for $N=2k$ even,  where $I$ is the left ideal generated by
the $\p_i$'s. $S$ has a basis $\x(m)$, where $m\in \R^k$ with $m_i\in\{ 0,1\}$
for $1\leq i\leq k$, and where
$$\x(m)=(\pd_1)^{m_1}(\pd_2)^{m_2}\ ...\ (\pd_k)^{m_k}\quad mod \ I.$$
If we assign to the highest weight vector $m(0)$ the weight $\epsilon =(\frac{1}{2}, \frac{1}{2},\ ...\ \frac{1}{2})$,
the vector $\x(m)$ would have the weight $\mu$ with $\mu_i=\frac{1}{2}-m_i$.
We also define the quantity
$$m\{ r\}=\sum_{j=1}^r m_j.$$
Then it follows that $\pd_j\x(m)=0$ if $m_j=1$, and 
$$\pd_j\x(m)=(-1)^{m\{ j-1\} }\x(\bar m^j), \quad {\rm if}\  m_j=0,$$
where $\bar m^j$ coincides with $m$ except in the $j$-th coordinate, which is
replaced by $1-m_j$. Observe that  the weight of $\x(\bar m^j)$ differs
from the one for $\x(m)$ only
in the $j$-th coordinate, by a sign. The action of $\p_j$ is given by the adjoint of $\pd_j$.

If $N=2k+1$ is odd, we can make $S$ as before into a $Cl(2k+1)$-module as follows.
The actions of $e_i$ with $i\leq 2k$  resp $\p_j$ and $\pd_j$ with $j\leq k$
is as before. The action of $e_{2k+1}$
is given by 
$$e_{2k+1}\x(m)\ =\ \pm f_{2k}\x(m)\ =\ \pm (-1)^{m\{k\}}\x(m),$$
where we obtain a $Cl(2k+1)$ action for each choice of the sign.

\subsection{Tensor products}\label{tensorprod} Let $S$ be the $Cl(N)$-module as described
in the previous subsection, for both $N$ odd and $N$ even.
In the following, we will identify elements of $Cl(N)$ with its image in $\End(S)$.
 We define
 elements $C \in \End(S^{\otimes 2})$ by
$$C=\frac{1}{2}\sum_{i=1}^N e_i\otimes e_i.$$
Using the definitions in \ref{psidef}, we  obtain
$$
\frac{1}{2}( e_{2j-1}\otimes e_{2j-1} + e_{2j}\otimes e_{2j})\ =\ \p_j\otimes \pd_j + \pd_j\otimes \p_j.$$
Using this and the definitions in the last section, we can also write for $N=2k+1$
\begin{equation}\label{Cclassic}
C\ =\ \frac{1}{2} f_{2k}\otimes f_{2k}\ +\ \sum_{j=1}^k \ \p_j\otimes \pd_j + \pd_j\otimes \p_j,
\end{equation}
while for $N=2k$ even $C$ is as above without the first summand. We will also need the well-known
isomorphisms of $Pin(N)$ modules given by 
\begin{equation}\label{Stensor}
N=2k:\quad S^{\otimes 2}\ \cong\ \bigoplus_{j=0}^{2k}\ \wedge^j V,\hskip 2em
N=2k+1:\quad S^{\otimes 2}\ \cong\ \bigoplus_{j=0}^{k}\ \wedge^j V,
\end{equation}
where $V=\C^N$ is the vector representation of $O(N)$.

\begin{lemma}\label{eigenvalueslemma} (a) If $N=2k$, the element $C$ has exactly the integer eigenvalues
$j$ satisfying $-k\leq j\leq k$.

(b) If $N=2k+1$, the element $C$ 
has the eigenvalues $\ve(-1)^j(k+\frac{1}{2}-j)$, $1\leq j\leq k$ with the sign $\ve$ depending on
the choice of the $Cl(N)$-module $ S$.
\end{lemma}

$Proof.$ Observe that 
$$ (\p_j\otimes \pd_j+ \pd_j\otimes \p_j) \x(m)\otimes \x(n)
\ =\ 
(-1)^{m\{ j-1\} -n\{ j-1\} }\ \x(\bar m^j)\otimes  \x(\bar n^j), $$
if $m_j+n_j=1$,
and it is equal to 0 otherwise. Let now $\rho=(k-i)_i\in\R^k$, and let $\bar m$ be defined by
$\bar m_j=1-m_j$, $1\leq j\leq k$. Let
\begin{equation}\label{eigenvector}
v\ =\ \sum_m (-1)^{(m,\rho)} \x(m)\otimes \x(\bar m).
\end{equation}
Then the $\x(m)\otimes \x(\bar m)$ coordinate of $Cv$ is given by 
$$(-1)^{(m,\rho)}\sum_j (-1)^{2m\{j-1\} -(j-1)+(\bar m^j-m,\rho)}= (-1)^{(m,\rho)}(-1)^{k-1}k.$$
This shows that $(-1)^{k-1}k$ is an eigenvalue for $C$. One obtains an eigenvector with eigenvalue
$(-1)^kk$ by multiplying the $\x(m)\otimes \x(\bar m)$-coordinate of the vector $v$ by
the scalar $(-1)^{m\{ k\} }$. One can similarly define eigenvectors $v$ with eigenvalues $\pm j$
by the sum of vectors $\x(m)\otimes \x(n)$ where $m_i+n_i=1$ for $1\leq i\leq j$ and where
$m(i)=n(i)=0$ for $i>j$.  These are all possible eigenvalues in view of \ref{Stensor} and 
Lemma \ref{Crelations}.

It follows from Lemma \ref{fnlemma},(c) that $f_N\x(m)\otimes f_N\x(n)=(-1)^{m\{ k\} + n\{ n\}}\x(m)\otimes \x(n)$.
It follows from this and statement (a) that the eigenvalues of $C$ are given by $(-1)^j(k+\frac{1}{2}-j)$, $1\leq j\leq k$.

\begin{lemma}\label{Crelations}  The element $C$   commutes
with the action of $Pin(N)$ on $S^{\otimes 2}$. Moreover, if 
 $C_1=C\otimes 1\in\End(S^{\otimes 3})$
and $C_2=1\otimes C\in\End(S^{\otimes 3})$ then we have
$$C_1^2C_2+2C_1C_2C_1+C_2C_1^2=C_2, \hskip 3em C_2^2C_1+2C_2C_1C_2+C_1C_2^2=C_1.$$
\end{lemma}

$Proof.$ As $g.(\sum_{i=1}^N e_i\otimes e_i)=\sum_{i=1}^N e_i\otimes e_i \in V^{\otimes 2}$ for any $g\in O(N)$,
the corresponding element $\sum_{i=1}^N e_i\otimes e_i\in Cl(n)^{\otimes 2}$ commutes with the action of $Pin(N)$ on
$ S^{\otimes 2}$. The first statement  follows from this and Lemma \ref{fnlemma},(d).

 Let $[A,B]_+=AB+BA$ for any two elements $A,B$ of a ring. Then it follows
from the relations that
\begin{equation}\label{eiejsum}
\sum_{i,j=1}^N\  [e_i,e_j]_+\ =\ \sum_{i,j=1}^N\ 2\delta_{i,j}1.
\end{equation}
Using this, we obtain
$$
[C_1,C_2]_+\ = \ \frac{1}{4}\sum_{i,j=1}^N\ e_i\otimes [e_i,e_j]_+\otimes e_j\ =\ 
       \frac{1}{2}\sum_{i=1}^N\ e_i\otimes 1 \otimes e_i.$$
We obtain in the same fashion
$$
[C_1,[C_1,C_2]_+]_+\ =\ \frac{1}{4}\sum_{i,j=1}^N\ [e_i,e_j]_+\otimes e_i\otimes e_j\ 
        =\ \frac{1}{4}\sum_{i,j=1}^r\ 2\delta_{ij}1\otimes e_i\otimes e_j\  =\ C_2.
$$
This proves the first identity in the statement. The proof of the second identity goes
the same way.

\section{ Representations of $U'_q\so_n$}\label{classicalcase}
We assume throughout this paper all the algebras to be defined over the field
of complex numbers, with $q$ not being a root of unity. See Section \ref{genring}
for more general rings.

\subsection{The algebras $U'_q\so_n$ and $U'_q\o_n$} We shall need a $q$-deformation
$U'_q\so_n$ of the universal enveloping algebra of the Lie algebra $\so_n$.
It is defined via generators $B_i$, $1\leq i<n$ and relations
$B_iB_j=B_jB_i$ for $|i-j|>1$ and
$$B_i^2B_{i\pm 1}-(q+q^{-1})B_iB_{i\pm 1}B_i+B_{i\pm 1}B_i^2=B_{i\pm 1},$$
with the choice of sign in the indices the same for all terms.
This algebra was defined, independently from each
other, by Gavrilik and Klimyk \cite{GK}, by Letzter \cite{Le} and by Noumi and
Sugitani \cite{NS}.  Its finite-dimensional representations
were classified by Klimyk and collaborators (see \cite{IK} and references there).

\begin{theorem}\label{Urepresentations} Let $q$ not be a root of unity.
Then there are two series of finite-dimensional
irreducible representations of $U'_q\so_n$, where $k=\lfloor n/2\rfloor$:

(a) The {\rm classical representations} are $q$-deformations of
the representations of $\so_n$. They are labeled by the dominant integral
weights of $\so_n$. They are given by all vectors
$\la=(\la_1,\la_2,\ ..., \la_k)$, where  all coefficients are either integers
or they all are congruent to $\frac{1}{2}$ mod $\Z$, and such that
$\la_1\geq \la_2\geq ...\ \geq\la_k\geq 0$ (for $n$ odd), or 
$\la_1\geq \la_2\geq ...\ \geq |\la_k|$ (for $n$ even). They have the same dimensions as the
corresponding $\so_n$ representations.

(b) The {\rm nonclassical representations} are labelled by all dominant
integral weights of $\so_n$ whose coefficients are not integers.
Their dimensions are $2^{-k}$ times the dimension of the corresponding
classical representations for $n$ odd, and $2^{1-k}$ times the dimension
of the corresponding classical representations for $n$ even. 
For each such weight, we have $2^{n-1}$ non-equivalent
representations of $U'_q\so_n$. They can be obtained from each other by
multiplying the matrices for various generators by $-1$. For $n$ even, it suffices to consider
only representations with highest weights $\la$ for which $\la_k>0$, together
with the just mentioned operation of sign changes.
\end{theorem}

\begin{remark}\label{GZremark} The construction of the representations 
of $U'_q\so_n$ by Klimyk et al essentially is a $q$-version of the construction
of representations of $\so_n$ via Gelfand-Zetlin bases, see \cite{Mo}.
This approach takes advantage of the fact that a simple $\so_n$-module, viewed
as an $\so_{n-1}$ module decomposes into a direct sum of mutually nonisomorphic
simple $\so_{n-1}$ modules, see e.g. \cite{Mo}, Section 4.1 and 4.2 for details.
This also determines the decomposition of a simple $classical$ $U'_q\so_n$ module
into a direct sum of mutually nonisomorphic $U'_q\so_{n-1}$ modules.

The decomposition of a simple {\it non-classical} $U'_q\so_n$ module  $V_\la$ into a direct sum
$\bigoplus V_\mu$ of simple $U'_q\so_{n-1}$ modules can be described
similarly: If $n=2k$ is even, $\mu$ runs through all weights $\mu$ satisfying
\begin{equation}\label{restricteven}
\la_1\geq\mu_1\geq \la_2\geq\ ...\ \geq\mu_{k-1}\geq \la_k>0,
\end{equation}
where for $n=2k+1$ odd we have
\begin{equation}\label{restrictodd}
\la_1\geq\mu_1\geq \la_2\geq\ ...\ \geq\mu_{k-1}\geq \la_k>\mu_k>0.
\end{equation}
All quantities in these inequalities are half-integers, i.e. congruent to 1/2 mod $\Z$.
For classical representations, the restriction rules are almost the same.
We only need to replace $\la_k$ by $|\la_k|$ for $n=2k$,
and $\mu_k$ by $|\mu_k|$ for $n=2k+1$ in the inequalities above,
see \cite{Mo}.
\end{remark}

We will also need an analog of the full orthogonal group in this setting. 
We define the algebra $U'_q\o_n$ by adding an additional generator $F$
to the generators of $U'_q\so_n$ with the relations
$$F^2=1,\quad FB_1=-B_1F\quad {\rm and}\quad FB_i=B_iF\quad {\rm for}\ i>1.$$

\begin{remark} 1. It is well-known that the Lie algebra $\so_n$ can be defined via
generators $L_i=E_{i,i+1}-E_{i+1,i}$, $1\leq i<n$, with $E_{i,j}$ matrix units.
It is then easy to check that the maps $B_j\mapsto \sqrt{-1}L_j$, $1\leq j<n$
 and $F\mapsto diag(-1,1,1,\ ...\ 1)$
define representations of $U'_1\so_n=U\so_n$ and $U'_1\o_n$ respectively.

2. It is also clear from the representation in the first remark that any irreducible
representation of the group $O(N)$ defines a representation of $U'_1\o_n$
by viewing the image of $F$ as a group element, and identifying $U'_1\so_n$ with
the universal enveloping algebra of $\so_n$. We will see in this paper that
these representations also exist and remain irreducible for generic $q$.
\end{remark}

\subsection{Homomorphism onto Temperley-Lieb algebras} The Temperley-Lieb algebra
$TL_n$ is given by generators $e_i$, $1\leq i<n$ and relations $e_i^2=e_i$, $e_ie_j=e_je_i$ for 
$|i-j|>1$ and $e_ie_{i\pm 1}e_i=\frac{1}{q^2+q^{-2}}e_i$.
It is well-known that  $\End_{U_{q^2}\sl_2}(V^{\otimes n})$ is isomorphic to $TL_n$
in our parametrization for $V=\C^2$. The proof of the following proposition is a straightforward,
if moderately tedious calculation. It is a special case of our main results Theorem \ref{fft}
and \ref{sft}.

\begin{proposition}\label{Temperley} The map $B_i\mapsto \frac{1}{q+q^{-1}}-(q+q^{-1})e_i$,
$1\leq i<n$ defines an algebra homomorphism from $U'_{-q^2}\so_n$ onto $TL_n$,
which induces non-classical representations of  $U'_{-q^2}\so_n$.
\end{proposition}

\subsection{Representations of $U'_q\so_3$}  Similarly as for the case of the Lie algebra $\sl_2$,
it is easy to write down explicit irreducible  representations for $U'_q\so_3$.
The $(N+1)$-dimensional irreducible  $classical$ representation $V_N$ of $U'_q\so_3$
can be described as follows. 
We fix a basis $\{v_j, 0\leq j\leq N\}$ of weight vectors. 
Then the actions of $B_1$ and $B_2$ are given by
\begin{equation}\label{so3action}
B_1v_j=[N/2-j]v_j,\quad B_2v_j=v_{j+1}+\al_{j-1,j}v_{j-1},
\end{equation}
where 
$$\al_{j-1,j}=\frac{[N+1-j][j]}{(q^{N/2-j}+q^{j-N/2})(q^{N/2-j+1}+q^{j-N/2-1})}.$$
Similarly, for $N$ odd, we can describe an $(N+1)/2$-dimensional simple 
{\it non-classical}  representation of $U'_q\so_3$ with respect to a basis $\{ v_j, 0\leq j< (N-1)/2\}$
by
\begin{equation}\label{so3actionnon}
B_1v_j=[N/2-j]^+v_j,\quad B_2v_j=v_{j+1}+\al^+_{j-1,j}v_{j-1},\quad j<(N-1)/2,
\end{equation}
where 
\begin{equation}\label{alphadef}
\al^+_{j-1,j}=\frac{[N+1-j][j]}{(q^{N/2-j}-q^{j-N/2})(q^{N/2-j+1}-q^{j-N/2-1})}.
\end{equation}
If $j=(N-1)/2$, the action on $v_j$ by $B_1$ is as in \ref{so3actionnon}, while we have
$$B_2v_{(N-1)/2}= \pm \frac{[(N+1)/2]}{i(q^{1/2}-q^{-1/2})}v_{(N-1)/2}\ +\ \al^+_{(N-3)/2,(N-1)/2}v_{(N-3)/2},$$
where the representation with the minus sign in the formula above is equivalent to the
representation with the plus sign, after replacing $B_2$ with $-B_2$. We will consider representations
of the algebra $U'_{-q^2}\so_n$. If we choose $(-q^2)^{1/2}=iq$, the eigenvalues of
$B_1$ in the $(N+1)/2=k+1$-dimensional nonclassical representation will be 
$(-1)^{k-j}(q^{N-2j}-q^{2j-N})/(q^2-q^{-2})$, $0\leq j\leq k$.

\begin{lemma}\label{so3rep} (a) We can make $V_N$ into a $U'_{-q}\so_3$-module
by leaving the action of $B_2$ the same, and replacing $B_1$ by $\tilde B_1$
whose action on $V_N$ is given by $\tilde B_1v_j=(-1)^jB_1v_j$.

(b) If $N$ is even, the representation of $U'_{-q}\so_3$ in (a) is isomorphic to
its classical representation with highest weight $[N/2]$. If $N$ is odd,
it decomposes into the direct sum of two non-classical  irreducible representations
of $U'_{-q}\so_3$ with highest weight $[N/2]^+=i(q^{N/2}+q^{-N/2})/(q-q^{-1})$.
\end{lemma}

$Proof.$ Obviously the actions of $B_1^2$ and $\tilde B_1^2$ on $V_N$ are the same,
while one checks easily that $\tilde B_1B_2\tilde B_1v_j=-B_1B_2B_1v_j$
for all basis vectors of $V_N$. This implies (a), using Theorem \ref{Urepresentations}.

If $N$ is even, it is easy to check that $v_N$ is a highest weight vector
with weight $[N/2]_{-q}$ for the $U'_{-q}\so_3$ action on $V_N$ given via
$\tilde B_1$ and $B_2$ which remains irreducible. For $N$ odd, $\tilde B_1$ acts
via the same scalar on $v_j$ as on $v_{N-j}$. Hence we obtain 
two highest weight vectors (in the sense  of \cite{WUrep})
with highest weight $[N/2]_+$, where the $q$-number here is defined
for $-q$. It follows that $V_N$ decomposes into the direct sum of two
irreducible non-equivalent $U'_{-q}\so_3$-modules  (see e.g. \cite{WUrep}, Theorem 3.8
for details).

\subsection{Homomorphisms}\label{homomorphisms}
 Let $S$ be the spinor module as described in 
Section \ref{explicit}
and let $C$ and $\tilde C$ be as in Section \ref{tensorprod}. We define elements
$$C_ i=\ 1\otimes 1 \otimes\ ... \otimes C \otimes\ ...\ \otimes 1\in \End(S^{\otimes n}), \hskip 3em 1\leq i<n,$$

$$\tilde C_ i=\ 1\otimes 1 \otimes\ ... \otimes \tilde C \otimes\ ...\ \otimes 1\in \End(S^{\otimes n}), 
\hskip 3em 1\leq i<n,\quad i \ {\rm odd},$$
where the element $C$ on the right hand side acts on the $i$-th and $(i+1)$-st factor of $S^{\otimes n}$.
Moreover, we define $\tilde C_{2i}=C_{2i}$  for $1< 2i<n$.

\begin{proposition}\label{prophom} For both $N$ odd and even, 
the map $B_i\mapsto C_i$ defines a homomorphisms of $U'_{-1}\so_n$ into
$\End_{Pin(N)}(S^{\otimes n})$. For $N$ even, the map $B_i\mapsto \tilde C_i$ for $i$ odd
and $B_i\mapsto  C_i$ for $i$ even defines a homomorphism of $U\so_n$ into 
$\End_{Pin(N)}(S^{\otimes n})$.
\end{proposition}

$Proof.$ The first statement follows from Lemma \ref{Crelations} and the discussion before this proposition.
The second statement was already shown in \cite{Wspin}. It can also be deduced from the first statement
using Lemma \ref{so3rep}.

\section{Quantum groups}

\subsection{$q$-Clifford algebra}\label{qCliffsection}
 We follow the paper \cite{Hay} by Hayashi, with some minor modifications.
For a somewhat more conceptual approach to $q$-Clifford algebras, 
see \cite{DF} and Section \ref{DFconn}.
The $q$-Clifford algebra $Cl_q(2k)$ coincides with the ordinary Clifford algebra
in the sense that it is again generated by elements $\p_i$ and $\pd_i$, $1\leq i\leq k$
satisfying the relations for the usual Clifford algebra.
In particular, $\p_i\pd_i$ and $\pd_i\p_i$ are idempotents which annihilate each other
and add up to 1.  The dependency on $q$ will be reflected by additional elements
$\omi$ defined by
\begin{equation}\label{omidef}
\omi\ =\ \p_i\pd_i+q^{-1}\pd_i\p_i.
\end{equation}
Then it is clear that $\omi^{\pm 1}\omj^{\pm 1}=\omj^{\pm 1}\omi^{\pm 1}$
and that
\begin{equation}\label{omirel}
\omi\p_i =\p_ i=q\p_i\omi, \hskip 3em \pd_i\omi=\pd_i=q\omi\pd_i.
\end{equation}

\subsection{Homomorphisms into $q$-Clifford algebras} We are now defining maps from
the quantum groups $U_q\so_{2k}$ and $U_q\so_{2k+1}$, i.e. of Lie type $D_k$ and $B_{k}$
into $Cl_q(N)$. Here we use the definition of the quantum groups as in
\cite{Lu} or in \cite{jantzen}, Section 4.3. In particular, the inner product on the weight lattice
is normalized such that $(\alpha_i,\alpha_i)=2$ for every $short$ root.
The maps appeared before in \cite{Hay}. However, our normalizations are not always the same
as in that paper, so we give the explicit maps below as 
follows: For $1\leq i\leq k-1$ we define
\begin{equation}\label{qdef1}
K_i\mapsto \om_i^2\om_{i+1}^{-2}, \hskip 2em E_i\mapsto \p_i\pd_{i+1}, 
\hskip 2em F_i\mapsto \p_{i+1}\pd_i.
\end{equation}
For type $B_k$, we also define
\begin{equation}\label{qdefB}
K_k\mapsto q\om_k^2, \hskip 2em E_k\mapsto \p_kf_{2k},\hskip 2em F_k\mapsto f_{2k}\pd_k,
\end{equation}
while for type $D_k$ we define
\begin{equation}\label{qdefD}
K_k\mapsto q^2\om_{k-1}^2\om_k^2, \hskip 2em E_k\mapsto \p_{k-1}\p_k,
\hskip 2em F_k\mapsto \pd_{k-1}\pd_k,
\end{equation}

\begin{proposition}\label{qgrouprep} (see \cite{Hay})
(a) The assignments in \ref{qdef1} and \ref{qdefD} define a representation of $U_{q^2}\so_N$, $N=2k$ even
into $Cl_q(N)\cong \End(S)$.

(b) The assignments in \ref{qdef1} and \ref{qdefB} define a representation of $U_q\so_{2k+1}$
into  $Cl_q(2k)\cong \End(S)$.
\end{proposition}

$Proof.$ Statement (a) was proved in \cite{Hay}. The assignments in (b)  differ from the ones
in \cite{Hay} only by multiplying the images of $E_{k}$ and $F_{k}$ by $f_{2k}$ from the right and
from the left. It is not hard to check that this still satisfies the quantum group relationis,
as $f_{2k}^2=1$ and $f_{2k}$ commutes with the images of the lower indexed generators.

\begin{remark}\label{quantumremark}
1. Observe that we have defined a representation of $U_{q^2}\so_N$ for $N=2k$. This will make
it easier to deal with the odd- and even-dimensional cases at the same time.

2. One can check that for $N=2k+1$ odd the vector $\x(0)$ is a highest weight vector
with weight $\ve=\frac{1}{2} (1,1,\ ..., 1)$, and that the vector $\x(m)$ has
weight $\mu$ with $\mu_i=\frac{1}{2}-m_i$, $1\leq i\leq k$. If $N=2k$ even,
we also have the highest weight $\ve_-=\x(m)$ with $m_i=\delta_{i,k}$. 
The weight of $\x(m)$ can be determined as in the odd-dimensional case.
\end{remark}

\subsection{Commuting objects, Lie type $D$} We now define the $q$-deformations of the operators
$C$ of the previous section for quantum groups. As before, we define them as elements of 
$Cl_q(N)^{\otimes 2}$ acting on
$\End(S^{\otimes 2})$. For Lie type $D_k$, we define
\begin{equation}\label{comqD}
C\ =\ \sum_{i=1}^k \ \Om_{i-1}^{-1}\p_i\otimes\Om_{i-1}\pd_i\ +\ \Om_{i-1}^{-1}\pd_i\otimes\Om_{i-1}\p_i,
\end{equation}
where $\Om_r=\prod_{j=1}^r\om_j^2$, $1\leq r\leq k$.
 We leave it to the reader to check that
\begin{equation}\label{Caction}
(\Om_{j-1}^{-1}\otimes \Om_{j-1})(\p_j\otimes \pd_j+ \pd_j\otimes \p_j) \x(m)\otimes \x(n)
\ =\ 
(-q^2)^{m\{ j-1\} -n\{ j-1\} }\ \x(\bar m^j)\otimes  \x(\bar n^j),
\end{equation}
if $m_j+n_j=1$,
and it is equal to 0 otherwise.

\begin{lemma}\label{comlemD}
The operator $C$ defined in Eq \ref{comqD} commutes with the action of $U_{q^2}\so_{2k}$ on $S^{\otimes 2}$.
\end{lemma}

$Proof.$ Let us first do this for type $D_2$. We shall use the coproduct defined by
$$\Delta(E_i)=K_i^{1/2}\otimes E_i+E_i\otimes K_i^{-1/2};$$
it is well-known that this is equivalent to the coproduct defined in \cite{Lu} and \cite{jantzen}, using the automorphism
defined by $E_i\mapsto E_iK_i^{1/2}$, $F_i\mapsto K_i^{-1/2}F_i$, $K_i\mapsto K_i$, $1\leq i\leq k$.
 Using Def. \ref{qdef1}, we obtain
$$\Delta(E_1)=\om_1\om_2^{-1}\otimes \p_1\pd_2+\p_1\pd_2\otimes\om_1^{-1}\om_2.$$
Let $C_1=\p_1\otimes\pd_1+\om_1^{-2}\pd_2\otimes\om_1^2\p_2$. We then obtain
\begin{align}
[ \Delta(E_1), C_1] 
=&\  \om_1\om_2^{-1}\p_1\otimes\p_1\pd_2\pd_1 -\p_1\om_1\om_2^{-1}\otimes\pd_1\p_1\pd_2\cr
&\ +\p_1\pd_2\om_1^{-2}\p_2\otimes\om_1^{-1}\om_2\om_1^2\pd_2
- \om_1^{-2}\p_2\p_1\pd_2\otimes\om_1^2\pd_2\om_1^{-1}\om_2\cr
=&\ -\om_2^{-1}\p_1\otimes(\p_1\pd_1+q^{-1}\pd_1\p_1)\pd_2 + \p_1(q\pd_2\p_2+\p_2\pd_2)\otimes\om_1\pd_2\cr
=&\ -\om_2^{-1}\p_1\otimes\om_1\pd_2\ +\ \om_2^{-1}\p_1\otimes\om_1\pd_2\ =\ 0,
\end{align}
where we used the relations \ref{omidef} and \ref{omirel}
after the definition of the $q$-Clifford algebra. One similarly also shows
that the commutant of $\Delta(E_1)$ with $\pd_1\otimes\p_1+\om_1^{-2}\p_2\otimes\om_1^2\pd_2$ is equal to 0.
This shows that $[\Delta(E_1),C]=0$. The statement for $F_1$ can be shown by a similar calculation, or it can
be deduced by the following argument: The 
transpose map $^T$ induced by 
$$\p_i\mapsto \p_i^T=\pd_i,\hskip 3em\pd_i\mapsto (\pd_i)^T=\p_i$$
induces an algebra antiautomorphism of $Cl_q$ which induces on the image of $U_q\so_N$
the algebra anti-automorphism defined by
$$E_i\mapsto E_i^T=F_i,\quad F_i\mapsto F_i^T=E_i,\quad K_i\mapsto K_i^T=K_i,$$
which is compatible with the Hopf
algebra structure. As $C^T=C$, it follows 
$$[\Delta(F_i),C]=[\Delta(E_i^T),C^T]=-[\Delta(E_i),C]^T=0.$$
The commutation relation of $C$ with $\Delta(E_2)$ is shown by a similar calculation, from which follows
the claim for $F_2$ by the previous argument.

For the general case, one observes that $\Delta(E_i)$ trivially commutes with all summands of $C$ except
the ones indexed by $i$ and $i+1$. The proof that $\Delta(E_i)$ commutes with these remaining summands
is essentially the same as the one for the case $D_2$.

\subsection{Commuting objects, Lie type $B$}
For Lie type $B_{k}$,  we define the operator $C\in \End(S^{\otimes 2})$   by (compare with
Section \ref{tensorprod} for $q=1$)
\begin{equation}\label{comqB}
C\ =\  \frac{1}{[2]}(\Om_{k}^{-1} f_{2k}\otimes \Om_{k} f_{2k})
\ +\ \sum_{i=1}^{k} \ \Om_{i-1}^{-1}\p_i\otimes\Om_{i-1}\pd_i\ +\ \Om_{i-1}^{-1}\pd_i\otimes\Om_{i-1}\p_i.
\end{equation}

\begin{lemma}\label{comlemB}
The operator $C$ defined in Eq \ref{comqB} commutes with the action of $U_q\so_{2k+1}$ on $S^{\otimes 2}$.
\end{lemma}

$Proof.$ The fact that the images of generators labeled by $i<k$ commute with $C$ except
for the first summand  follows from the proof of Lemma \ref{comlemD}. One then checks that
$\p_i\pd_{i+1}$ commutes with $\om_i\om_{i+1}$ by a direct calculation; this implies that it
also commutes with $\Om_k$. It is now easy to check that $\Delta(E_i)$ also commutes
with the first summand in the definition of $C$.

To finish the proof, recall that we have (with $N=2k$)
$$\Delta(E_{k})\ =\  q^{1/2}\om_{k}\otimes \p_{k}f_N\ +
\ \p_{k}f_N\otimes  q^{-1/2}\om_{k}^{-1}.
$$
We will use relations \ref{omidef} and \ref{omirel}, which also imply
$\p_{k-1}\Om_{k-1}=q^{-2}\Om_{k-1}\p_{k-1}$. 
We obtain
\begin{align}
& [\Delta(E_{k}), \Om_{k}^{-1} f_N\otimes \Om_{k-1} f_N]\ = \cr
= & 
\ q^{1/2}\om_{k} \Om_{k}^{-1} f_N\otimes  [\p_{k}f_N\Om_{k} f_N-\Om_{k} f_N\p_{k}f_N] \cr
& +\ [\p_{k}f_N\Om_{k}^{-1}f_N- \Om_{k}^{-1}f_N\p_{k}f_N] \otimes   q^{-1/2}\om_{k}^{-1}\Om_{k} f_N\cr
= &\   \om_{k} \Om_{k}^{-1} f_N\otimes (q^{-3/2}+q^{1/2})\Om_{k}\p_{k}
+ (q^{3/2}+q^{-1/2})\Om_{k}^{-1}\p_{k}\otimes \Om_{k} f_N \om_{k}^{-1}\cr
= &\ (q+q^{-1})\  (q^{-1/2}\om_{k} \Om_{k}^{-1} f_N\otimes \Om_{k}\p_{k}
+ q^{1/2} \Om_{k}^{-1}\p_{k}\otimes \Om_{k-1} f_N \om_{k-1}^{-1}).
\end{align}
We also have
\begin{align}
& [\Delta(E_k), (\Om_{k-1}^{-1}\otimes \Om_{k-1})(\p_{k}\otimes\pd_{k} + \pd_{k}\otimes\p_{k})]\ =\cr
=& \ q^{1/2}\om_{k}\Om_{k-1}^{-1}\p_{k}\otimes \p_{k}f_N\Om_{k-1}\pd_{k}
- \Om_{k-1}^{-1}\p_{k} \ q^{1/2}\om_{k}\otimes \Om_{k-1}\pd_{k}\p_{k}f_N\cr
& + \p_{k}f_N\Om_{k-1}^{-1}\pd_{k}\otimes q^{-1/2}\om_{k}^{-1}\Om_{k-1}\p_{k}\ -\ \Om_{k-1}^{-1}\pd_{k} \p_{k}f_N\otimes  q^{-1/2}\Om_{k-1}\p_{k}\om_{k}^{-1}\cr
=&\ -\Om_{k-1}^{-1}\p_{k}\otimes \Om_{k-1}f_N(q^{1/2}\p_{k}\pd_{k}+q^{-1/2}\pd_{k}\p_{k})\cr
&-\ (q^{-1/2}\p_{k}\pd_{k}+q^{1/2}\pd_{k}\p_{k})\Om_{k-1}^{-1}f_N\otimes\Om_{k-1}\p_{k}\cr
=&\ -  q^{1/2}\Om_{k}^{-1}\p_{k}\otimes\Om_{k-1}f_N\om_{k} -
 q^{-1/2}\om_{k}^{-1}\Om^{-1}_{k-1}f_N\otimes\Om_{k-1}\p_{k}.
\end{align}
Obviously, $\Delta(E_{k})$ commutes with the first $k-1$ summands under the summation sign in the definition of $C$.
It follows from the identity $\om_k\Om_{k}^{-1}=\om_k^{-1}\Om_{k-1}^{-1}$ and the last two calculations
that $\Delta(E_k)$ also commutes with the remaining summands of $C$.
The commutation with $\Delta(F_{k})$ can be deduced from this using the transposition map $^T$ as
in the proof of Lemma \ref{comlemD}.

\section{Relations}\label{relationssec}

The main result of this section will be to give an algebraic description of the centralizer of
the action of $U_q\so_N$ on $S^{\otimes n}$. It is possible, and fairly straightforward,
to extend the proof in \cite{Wspin} to the additional cases treated here. However, that proof 
was somewhat indirect. So we decided to give another proof here which, basically, is
a direct calculation. While not quite as straightforward as the proof for $q=1$
in Lemma \ref{Crelations}, it would still seem to be an improvement over the one in \cite{Wspin}.

\subsection{Basic relations} Let $a,b$ be any elements in an associative algebra,
and let $v$ be any invertible element in its ground ring. Then we define
\begin{equation}\label{lhs}
lhs_v(a;b)\ =\  a^2b+(v+v^{-1})aba+ba^2.
\end{equation}
 It will be convenient to introduce the notation
$$c_{i,+}=\Omega^{-1}_{i-1}\p_i,\quad c_{i,-}=\Omega^{-1}_{i-1}\pd_i,\quad 
d_{i,+}=\Omega_{i-1}\p_i,\quad d_{i,-}=\Omega_{i-1}\pd_i.$$
Using this, we can write the commuting operators $C$ from Lemma \ref{comlemD} and Lemma \ref{comlemB}
as
\begin{equation}\label{Cidef}
C\ =\ \sum_{i=1}^k\ C(i)\quad {\rm and}\quad
C\ =\ \tilde C(k+1)\ +\ \sum_{i=1}^k\ C(i) 
\end{equation}
where
$$C(i)\ =\ c_{i,+}\otimes d_{i,-}\ +\ c_{i,-}\otimes d_{i,+},\hskip 3em
\tilde C(k+1)\ =\ \frac{1}{[2]}\ \Om_k^{-1}f_N\otimes\Om_kf_N.$$
It is straightforward to check the following relations, where $\ve,\kappa \in\{ \pm\}$ and $q^{2\ve}=q^{\pm 2}$:

\begin{equation}\label{cdrel}
d_{i,\ve}c_{j,\kappa}\ =\ 
\begin{cases} -q^{2\ve } c_{j,\kappa}d_{i,\ve} & i<j,\cr
 -q^{2\kappa } c_{j,\kappa}d_{i,\ve} & i>j.
\end{cases}
\end{equation}
We obtain from the relations so far the following equation which will be useful later:

$$d_{i,\ve}d_{i,-\ve}c_{j,\kappa} +(q^2+q^{-2})  d_{i,\ve}c_{j,\kappa}d_{i,-\ve}  + c_{j,\kappa}d_{i,\ve}d_{i,-\ve}\ =\ $$
\begin{equation}\label{iirelation}
=\ \begin{cases} 0  & i>j, \cr  (1-q^{4\ve}) d_{i,\ve}d_{i,-\ve}c_{j,\kappa} & i<j.
\end{cases}
\end{equation}

\subsection{Technical lemma}
\begin{lemma}\label{Cirelations} Using notation defined in \ref{lhs} we have

(a) $lhs_{q^2}(C\otimes 1; 1\otimes C)\ =\ \sum_{i,k} lhs_{q^2}(C(i)\otimes 1; 1\otimes C(k))$,

(b) 
$$lhs_{q^2}(C(i)\otimes 1;1\otimes C(k))\ =\ 
\begin{cases} 0& i>k\cr
(\Om_{i-1}^{-2}\otimes \Om_{i-1}^2\otimes 1)(1\otimes C(k))& i=k,\cr
([\Om_{i-1}^{-2}\otimes \Om_{i-1}^2-\Om_{i}^{-2}\otimes \Om_{i}^2]\otimes 1))(1\otimes C(k))& i<k.
\end{cases}
$$
\end{lemma}

$Proof.$ It will be convenient to write $C=\sum_u c_u\otimes d_u$ for (a). (This is not quite
consistent with our previous notation, but should not lead to confusion). 
It then follows that
\begin{equation}\label{expansion}
lhs_{q^2} (C\otimes 1; 1\otimes C)\ =\ \sum_{u,v,w} c_uc_v\otimes [d_ud_vc_w+(q^2+q^{-2})d_uc_wd_v+c_wd_ud_v]\otimes d_w.
\end{equation}
Let $c_u=c_{i,\ve}$ and let $c_v=c_{j,\tve}$ with $i\neq j$. Observe that the claim
is proved if for given index $w$ the summand for our given indices $u$ and $v$ cancels with the one
with $c_u=c_{j,\tve}$ and $c_v=c_{i,\ve}$. Let $d_u=d_{i,-\ve}$ and $d_v=d_{j,-\tve}$,
and $c_u$ and $c_v$ as at the beginning of this paragraph.
It follows from \ref{cdrel}  that 
$$c_uc_v=-q^{\pm 2\ve}c_vc_u, \quad d_ud_v=-q^{\pm 2\ve}d_vd_u,$$
with matching signs in the exponents.   Let us choose the labeling such that $c_uc_v=-q^2c_vc_u$,
and hence also $d_ud_v=-q^2d_vd_u$.
Using the commutation relations above these two summands add up to
$$c_uc_v\otimes (q^2+q^{-2})[q^{-2}d_ud_vc_w+d_uc_wd_v-q^{-2}d_vc_wd_u+q^{-2}c_wd_ud_v]\otimes d_w.$$ 
Then our claim will follow if we can show that the middle factor $M$ in this tensor product is equal to 0.

First observe that $d_ud_v=-q^2d_vd_u$ implies that $d_u=d_{i,+}$ and $d_v=d_{j,\pm}$
or $d_v=d_{i,-}$ and $d_u=d_{j,\pm}$ with $i<j$. Let us consider the case $d_u=d_{i,+}=\om_{i-1}\p_i$.
Then one checks that 
$c_wd_u=-q^{-2}d_uc_w$ is only possible for $c_w=c_{a,-}$  for $a<i$.
As $i<j$, we then also have $c_wd_v=-q^2d_vc_w$, which forces $M=0$.
One similarly shows in the second case with $d_v=d_{i,-}$ that $c_wd_v=-q^2d_vc_w$ would imply $c_wd_u=-q^2d_uc_w$.
This completes the proof for claim (a).

For part (b), we only need to consider the cases with $c_u=c_{i,\ve}$ and $c_v=c_{i,-\ve}$, as $c_{i,\ve}^2=0$.
Using \ref{cdrel}, one checks
that
$$c_{i,\ve}c_{i,-\ve}\otimes [ d_{i,-\ve}d_{i,\ve}c_{j,\kappa}+(q^2+q^{-2})d_{i,-\ve}c_{j,\kappa}d_{i,\ve}+ c_{j,\kappa}d_{i,-\ve}d_{i,\ve} ]\otimes d_{j,-\kappa}\ =$$
$$=\ \begin{cases} 0 &  i>j,\cr
(1-q^{4\ve}) (\Om_{i-1}^{-2}\p_i^{\ve}\p_i^{-\ve}\otimes \Om_{i-1}^2\p_i^{-\ve}\p_i^{\ve}\otimes 1)(1\otimes c_{j,\kappa}\otimes d_{j,-\kappa})& i<j.
\end{cases}
$$
Adding up these quantities for all possible choices of $\ve$ and $\kappa$, we obtain 0 for $i>j$. Using
$$(1-q^4)\p_i\pd_i\otimes \pd_i\p_i+(1-q^{-4})\pd_i\p_i\otimes \p_i\pd_i= 1\otimes 1-\om_i^4\otimes\om_i^{-4},$$
we similarly obtain the claim for $i<j$. The claim for $i=j$ follows from a direct calculation.

\begin{proposition}\label{relationprop}
Let $C_1=C\otimes 1$ and $C_2=1\otimes C$. Then we have
$$C_1^2C_2+(q^2+q^{-2})C_1C_2C_1+C_2C_1^2= C_2.$$
\end{proposition} 

$Proof.$ It follows from Lemma \ref{Cirelations}(b) that for fixed $j$
$$\sum_{i=1}^k lhs_{q^2} (C(i)\otimes 1;1\otimes C(j))\ =\ 1\otimes C(j).$$
The claim follows from this and Lemma \ref{Cirelations}(a) for type $D$. 
For Lie type $B_k$, we have to add $\tilde C(k+1)$ to the expression
for type $D_k$, see \ref{Cidef}.
Setting $c_{k+1,\ve}=\Om^{-1}_{k}f_N$ and $d_{k+1,\ve}=\Om_k f_N$,
one checks that Eq \ref{cdrel} and \ref{iirelation} also hold if one
of the indices is $k+1$. 
One deduces the results of Lemma \ref{Cirelations}(a) and, except for $i=j=k+1$,
also of part (b) from this. Observe that $c_{k+1,\ve}$ commutes with $d_{k+1,\ve}$.
One calculates that
$$lhs_{q^2}(\tilde C(k+1)\otimes 1; 1\otimes \tilde C(k+1))= (\Om_k^{-2}\otimes \Om_k^2\otimes 1)(1\otimes \tilde C(k+1)).$$
The claim can be deduced from this, using Lemma \ref{Cirelations}, as it was done for type $D$.

\section{First and Second Fundamental Theorem}

\subsection{Preliminaries}  We consider the spinor module $S$ 
as in Section \ref{spinreps} for $\U=U_q\so_{2k+1}$,
for  $\U=U_q\so_{2k}$ and for $\U=U_q\so_{2k}\rtimes Z/2$,
where the $\Z/2$ action is given by the diagram automorphism
permuting the generators given by the end vertices of the Dynkin diagram
next to the triple vertex. It has dimension $2^k$
in all these cases, and its weights are given by all possible vectors
$\om =(\pm \frac{1}{2}, \pm \frac{1}{2},\ ...,\ \pm \frac{1}{2})\in \R^k$.
As all weights have multiplicity 1, tensoring an irreducible highest weight
module $V_\la$ by $S$ is given by
\begin{equation}\label{tensorrule}
V_\la\otimes S\ \cong\ \oplus_\mu V_\mu,
\end{equation}
where the summation goes over all dominant weights $\mu$ of the
form $\mu=\la + \om$.  If $\U=U_q\so_{2k}\rtimes Z/2$,
this has to be slightly modified, see \cite{Wspin}.
A first fundamental theorem has been proved for $\End_\U(S^{\otimes n})$ 
 for this case as well as for $\U=U_q\so_{2k+1}$
 in \cite{Wspin}, Theorem 3.3.
We will review the method used there, a modification of which will be used
also for the missing case to be proved here.
In all of these cases, we have 
\begin{equation}\label{tensordecomp}
S^{\otimes n}\ =\ S^{\otimes n}_{old}\oplus\ S^{\otimes n}_{new},
\end{equation}
where $S^{\otimes n}_{old}$ is a direct sum of irreducible
modules $V_{\lambda}$ which have already appeared in smaller tensor powers of $S$ 
(for which $\la_1< n/2$) and where $S^{\otimes n}_{new}$ is a direct sum
of irreducible modules $V_\la$ which have
not appeared before (for which $\la_1=n/2$).
The following proposition is a slight generalization of 
\cite{exc}, Prop. 4.10. We assume $\Ca$ to be a
$\C$-linear rigid braided tensor category,
see e.g. \cite{Turaev} for precise definitions. The only
property we will need is the fact that
for every object $X$ in $\Ca$ there exists an object $\bX$
and morphisms $\iota_X : \1\to X\otimes \bX$, $\tilde d_X: X\otimes \bX\to \1$
such that $Tr(a)= \tilde d_X (a\otimes 1) \iota_X$ for any 
$a\in \End(X)$. It is well-known that this holds for $\Ca =$ Rep$U_q\so_N$.

\begin{proposition}\label{basiccon}
Let $V$ be a self-dual rigid object in the braided spherical tensor category $\Ca$
and let $p=\iota_V\tilde d_V$. Let $\E_n=\End_\Ca(V^{\otimes n})$.
Then $\End(V^{\otimes n})_{old}= (\E_{n-1}\otimes 1)p_{n-1}(\E_{n-1}\otimes 1)$.
\end{proposition}

 $Proof.$ The $p$ in the statement can be normalized to be a projection.
It then satisfies exactly the same properties as the $p$ in the proof of
Proposition 4.10 in \cite{exc} for $k=2$, even if $V$ is not necessarily simple.
The claim follows from this.

\subsection{First fundamental theorem} We study $S^{\otimes n}_{new}$
by induction on the rank of $\U$.
It follows from the relations that the subalgebra of $U_q\so_N$,
generated by $E_i$, $F_i$ and $K_i^{\pm 1}$, $2\leq i\leq k$
is isomorphic to $U_q\so_{N-2}$, for $N=2k$ or $N=2k+1$.
It will be convenient to denote $U_q\so_N$ or $U_q\so_{N}\rtimes \Z/2$
by $\U(N)$,  the spinor module $S$  of $\U(N)$
by $S(N)$, and the centralizer algebra $\End_\U(S^{\otimes n})$
by $\E_n^{(N)}$. 
We have the following well-known facts which are easy to check:

(a) We have the isomorphism 
of $\U(N-2)$ modules $S(N)\cong S(N-2)_1\oplus S(N-2)_2$,
where $S(N-2)_1$ is spanned by the weight vectors with weights
$(\frac{1}{2},\om')$, with $\om'$ a weight of $S(N-2)$.

(b) Let $v\in S(N-2)_1^{\otimes n}\subset S(N)^{\otimes n}$.
Then $v$ is a highest weight vector for $\U(N-2)$ of weight
$\la'$  if and only if
it is a highest weight vector for $\U(N)$ of weight $(\frac{n}{2},\la')$.

\begin{theorem}\label{fft} (First Fundamental Theorem)
The algebra $\End_\U(S^{\otimes n})$ is generated by the elements
$a_i$, $i=1,2,\ ...\  n-1$, with $a\in \End_\U(S^{\otimes 2})$.
\end{theorem}

$Proof.$ We only need to consider the case $\U=U_q\so_{2k}$. 
The statement was proved in \cite{Wspin}, Theorem 3.3
for the other cases.

For getting the induction on the rank $k$ going, we define $Spin(2)$ to be the $\Z/2$
cover of $SO(2)$. Its irreducible representations are labeled by half
integers. In this case $S$ is the direct sum of two 1-dimensional representations
with weights $\pm \frac{1}{2}$, which is obviously self-dual.
Let $f\in \End(S)$ act via $\pm 1$ on the vectors with weights
$\pm \frac{1}{2}$. It follows from the tensor product rules
that $f\otimes 1_{n-1}$ again acts via $\pm 1$ on the 1-dimensional
representations with weights $\pm \frac{n}{2}$.
The tensor product rules also show that this is $S^{\otimes n}_{new}$.
The claim now follows for $SO(2)$ from this and Proposition
\ref{basiccon} by induction on $n$.

We similarly prove the claim for $N=2k>2$ by induction on $n$.
For $n=1$, $ S$ is the direct sum of two irreducible modules $S_\pm$,
on which the endomorphism $f$ acts via $ \pm 1$.
The claim follows  for $\End_\U(S^{\otimes n})_{old}$ from Proposition \ref{basiccon}
by induction assumption on $n-1$. 

It follows from observation (b) before this theorem that we have a surjective map
from $\End_{\U(N)}(S^{\otimes n})$ onto $\End_{\U(N-2)}(S_1^{\otimes n})$,
given by restriction from $S^{\otimes n}$ to $S_1^{\otimes n}$.
Indeed, the simple module of  $\End_{\U(N)}(S^{\otimes n})$
consisting of highest weight vectors of weight $(\frac{n}{2},\la')$
coincides with the simple module of $\End_{\U(N-2)}(S_1^{\otimes n})$
consisting of highest weight vectors of weight $\la'$.
Its kernel is $\End_{\U(N)}(S^{\otimes n}_{old})$.
By induction assumption, $\End_{\U(N-2)}(S_1^{\otimes n})
\cong \End_{\U(N)}(S^{\otimes n}_{new})$ is generated
by $\End_{\U(N-2)}(S_1^{\otimes 2})$. But the latter is just the
restriction of  $\End_{\U(N)}(S^{\otimes 2})$ to $S_1^{\otimes 2}$.
Hence   $\End_{\U(N)}(S^{\otimes 2})$ also generates $\End_{\U(N)}(S^{\otimes n}_{new})$.
The proof for $N$ odd goes the same way, with the case $N=3$ already
proved in Proposition \ref{Temperley}.

\subsection{Second fundamental theorem} Recall that the algebra
$U'_q\o_n$ was defined by adding an additional generator $F$
to $U'_q\so_n$ with the relations $F^2=1$, $FB_1=-B_1F$ and $FB_i=B_iF$ for $i>1$.
Also recall that the finite dimensional representations of the group $O(n)$
are labeled by all Young diagrams whose first two columns contain at most
$n$ boxes, see \cite{Wy} or also e.g. \cite{Wspin}. 

\begin{theorem}\label{sft} (Second Fundamental Theorem)
(a) If $N$ is odd and $\U= U_q\so_N$, we have a surjective map $U'_{-q^2}\so_n\to \End_{\U}(S^{\otimes n})$
defined by $B_i\mapsto C_i$,
where $C_i$ is defined as in Section \ref{homomorphisms},
using the map $C$ defined in \ref{comqB}. Its image is the direct sum of all non-classical representations
of  $U'_{-q^2}\so_n$  with highest weights $\mu$ such that $\mu_1\leq N/2$ 
and in  which all
$B_i$s have eigenvalues contained in $\{(-1)^j (q^{N-2j}-q^{2j-N})/(q^2-q^{-2}),\ 0\leq j<N/2\}$.

(b) If $N$ is even and $\U= U_q\so_N\rtimes Z/2$, we have a surjective map $U'_{-q}\so_n\to \End_{\U}(S^{\otimes n})$
defined by $B_i\mapsto C_i$, using the map $C$ defined in \ref{comqD}. 
Its image is the direct sum of all classical representations
of  $U'_q\so_n$  with highest weights $\mu$ such that $\mu_1\leq N/2$, with $\mu_i\in\Z$.

(c)  If $N$ is even and $\U= U_q\so_N$, we have a surjective map $U'_{-q}\o_n\to \End_{\U}(S^{\otimes n})$
defined by $B_i\mapsto C_i$ and by $F\mapsto f\otimes 1_{n-1}$. Its image is the direct sum of irreducible representations
of  $U'_q\o_n$  which specialize to representations of $O(N)$ labeled by Young diagrams $\mu$ whose first two columns contain
$\leq n$ boxes and such that $\mu_1\leq N/2$ for $q=1$.
\end{theorem}

$Proof.$  It was shown  in Proposition \ref{relationprop} that $B_i\mapsto C_i$ 
defines a representation of $U'_{-q^2}\so_n$.  Substituting
$q^2$ by $q$ for $N$ even, we get the representations as
stated. The eigenvalues of $C$ were computed
in Lemma \ref{eigenvalueslemma} for $q=1$. For general $q$,
the eigenvalues can be computed by the same method (see \cite{Wspin},
Lemma 4.2 and Proposition 4.3) or by a deformation argument,
using the explicit classification of irreducible representations of
$U'_q\so_3$.  From this also follows that 
$C$, respectively $C$ and $f\otimes 1$ for $\U=U_q\so_{2k}$
generate $\End_\U(S^{\otimes 2})$. The surjectivity statement now
follows from Theorem \ref{fft}. The explicit combinatorics will be studied
in the next section.

\subsection{Combinatorics and representations of $U'_q\so_n$ and $U'_q\o_n$} 
Our duality result implies that we can associate to each irreducible representation 
$V_\la$ of,
say, $\U=U_q\so_N\rtimes \Z/2$ which appears in $S^{\otimes n}$ 
an irreducible representation $W_{\la^c}$ of $U'_q\so_n$ such that
the multiplicity of $V_\la$ in $S^{\otimes n}$ is given by the dimension of $W_{\la^c}$.
The fact that such multiplicities are given by  dimension 
formulas (modified for $N$ odd) has been known for a long time, see e.g. \cite{Ba}. 
The Young diagram (or weight) $\la^c$ can be
obtained from the Young diagram $\la$ as its complement
in a rectangle with side lengths $N/2$ and $n/2$,
reflected at the line $y=x$, to get it into usual Young diagram
position. As our weights also
involve half integers, and there are some additional subtleties for $N$ even,
we will spell this out in more detail in this section, even though it is
not new (see e.g. \cite{Ba}, \cite{Has} or \cite{Wspin}).

Let us briefly describe the irreducible representations of $Pin(N)$.
If $N$ is even, the irreducible representations of $O(N)$ are labeled by Young diagrams,
whose first two columns contain at most $N$ boxes. The irreducible representations
of $Pin(N)$ which do not factor over $O(N)$ are given by $N/2$-tuples
$(\la_i)$ with $\la_1\geq \la_2\geq\ ...\ \geq \la_{N/2}>0$
with all $\la_i\equiv 1/2$ mod $\Z$. We also view them as Young diagrams
with an additional half box of width 1/2 and height 1 at the end.
We can now give a precse description of the diagram $\la^c$ which
has already appeared before, e.g. in \cite{Ba}.

\begin{definition}\label{complement} Let $R$ be a rectangle of height $N/2$ and width $n/2$.

(a)  Let $N=2k$ and let $\la$ be a Young diagram labeling an
irreducible representation of $Pin(N)$ appearing in $S^{\otimes n}$.
Then its $n$-complement  is the dominant integral weight
$\la^c$ of $\so_n$ such that $\la^c_j$ is equal to the number of boxes
in the $j$-th coulmn from the right in $R\backslash \la$.
If $\la$ does not fit into $R$, i.e. if the first column of $\la$ contains
more than $N/2$ boxes, we define $\la^c_{n/2}=n/2-\la'_1<0$,
where $\la'_1$ is the number of boxes in the first column of $\la$.

(b) If $N=2k$ is even and $\la$ a dominant integral weight labeling
an irreducible $\so_N$ module appearing in $S^{\otimes n}$,
we associate to it the Young diagram $\la^c$ whose $j$-th column
contains $n/2-\la_{N/2+1-j}$ boxes.

(c) For $N$ odd, and $\la$ an integral dominant weight labeling a representation
of $so_N$ appearing in $S^{\otimes n}$, we define the integral dominant
weight $\la^c$ labeling a non-classical representation of $U'_q\so_n$
by $\la^c_j$ to be equal to the number of boxes in the $j$-th column
from the right of $R\backslash \la$.
\end{definition}

\begin{proposition}\label{explduality} (see e.g. \cite{Ba}, \cite{Has}) Let $S$, $\U$, $U'_q\so_n$ etc be as at the
beginning of this section. Then the module $S^{\otimes n}$ has a multiplicity
one decomposition
$$S^{\otimes n}\ \cong\ \bigoplus_\la V_\la\otimes W_{\la^c},$$
where $\la$ ranges over the equivalence classes of irreducible $\U$
modules $V_\la$ which appear in $S^{\otimes n}$, and $\la^c$ is the 
label of the simple $U'_q\so_n$
resp $U'_q\o_n$-module $W_{\la^c}$ associated to it in Def \ref{complement}.
\end{proposition}

$Proof.$ We give a proof for the case $N$ odd by induction on $n$. 
We define $\la^r=\la^{r(n)}$ by $\la^r_i=n/2-\la_i$. This describes the 
number of boxes in the $i$-th row of the (unreflected) complement of $\la$
in $R$. 
Let $W_{\la^d}$ be the $U'_q\so_n$ module associated to $\la$
(e.g. we could take the module of highest weight vectors in $S^{\otimes n}$
with weight $\la$). It follows from the tensor product rules and induction assumption
that 
$$
W_{\la^d}\ \cong\ \bigoplus W_{\mu^{c(n-1)}}
$$
as $U'_q\so_{n-1}$-modules, where the summation goes over all highest weights
$\mu$ with $V_\mu\subset S^{\otimes n-1}$ such that $V_\la\subset V_\mu\otimes S$. 
This implies that
$\la=\mu+\om$ for some weight $\om$ of $S$. As $\om=\ve-\sum_j \ep_{i_j}$
for some subset $\{ i_j\}\subset \{ 1,2,\ ...,\ k\}$,
it follows that 
$$
\mu^{r(n-1)}= (n-1)\ve-\mu = n\ve-\la -\sum_j \ep_{i_j} = \la^{r(n)}-\sum_j \ep_{i_j},
$$
i.e. the rows of $\mu^{r(n-1)}$ and $\la^{r(n)}$ differ by at most 1.
Recall that $\la^{c(n)}$ is obtained from $\la^{r(n)}$ by reflecting the latter at 
the line $y=x$,
i.e. its $j$-th column contains $\la^{r(n)}_{(N+1)/2-j}$ boxes.
It is probably easiest seen geometrically that the condition
$0\leq \la^r_i-\mu^r_i\leq 1$, $1\leq i\leq N$, is equivalent to the branching conditions
\ref{restricteven} and \ref{restrictodd} for $\mu^c$ and $\la^c$. Hence the simple $U'_q\so_n$ module
$W_{\la^d}$ coincides with the simple  $U'_q\so_n$ module $W_{\la^c}$ as
a $U'_q\so_{n-1}$-module. It is easy to see from \ref{restricteven} and \ref{restrictodd}
that non-isomorphic simple $U'_qso_n$-modules are still non-isomorphic if viewed
as $U'_q\so_{n-1}$ modules. Hence $W_{\la^d}\cong W_{\la^c}$.

The other cases can be shown similarly. Another method would be to calculate
the character of $O(N)\times SO(n)\subset O(Nn)$ on the spin module of
$O(Nn)$. This was essentially done in \cite{Has}, with some minor errors.
Hopefully the statement of the pairings is correct in this paper.

\begin{remark}\label{symplectic} It was already observed in \cite{Ba}
that for $N$ odd the multiplicity of $V_\la$ in $S^{\otimes 2k+1}$ coincides with
the dimension of the irreducible $Sp(2k)$ module $W_{\la^c-\ve}$,
where $\ve=(\frac{1}{2},\frac{1}{2},\ ...,\ \frac{1}{2})\in\R^k$. No such
classical interpretation seems to be available for multiplicities in even tensor powers
of $S$.
\end{remark}

\begin{corollary}\label{onrep} Let $q$ not be a root of unity. 
Then every finite-dimensional irreducible representation of
$U'_q\so_n$  appears
in some tensor power $S^{\otimes n}$ for some $\U=U_q\so_N$
or $\U=U_q\so_N\rtimes Z/2$. Moreover,
every irreducible  finite dimensional  representation of $O(n)$ extends to an
irreducible representation of  $U'_q\o_n$.
\end{corollary}

$Proof$ This was shown for classical representations with integer highest
weights  in \cite{Wspin} and also in this paper. It can be shown the same way
for classical representations with half integer highest weights for $N$ odd, using the 
bigger module $\tilde S$ and \cite{Wspin} Lemma  4.2 and Proposition 4.3.
Existence of non-classical representations and of representations
of $U'_q\o_n$ as stated follows from Proposition
\ref{explduality}.

\subsection{Results for more general rings}\label{genring}  It is well-known that the Drinfeld-Jimbo quantum
groups can be defined over the ring $\Z[q,q^{-1}]$, see \cite{Lu}.
It is not hard to check that the same is true for the subalgebra
$U'_q\so_n\subset U_q\sl_n$. Also observe that the representations of $U'_q\so_n$
in Theorem \ref{sft} can be defined over the ring $\Z[q,q^{-1}]$
in type $D$ (for $N$ even) and over the ring $\Z[q,q^{-1}, [2]^{-1}]$ for $N$ odd.
Indeed, one checks easily that the matrix coefficients of the maps $C$ in 
\ref{comqD} and \ref{comqB} with respect to the basis $\x(m)$ are in
the rings as claimed.
We can deduce the following result from this and Theorem \ref{sft}:

\begin{proposition}\label{generalring}
The representations of the algebras in Theorem \ref{sft} on $S^{\otimes n}$
are well-defined also over the ring $\Z[q,q^{-1}]$ for $N$ even.
For $N$ odd, they are well-defined over the ring $\Z[q,q^{-1},[2]^{-1}]$.
The statements of Theorem \ref{sft} also hold in these cases, except
possibly the surjectivity statements.
\end{proposition}

\subsection{Results for fusion categories} For $q\neq \pm 1$ a root of unity 
the representation category of a quantum group $U_q\g$ is not semisimple.
But one can define an important semisimple
quotient tensor category $\F$  of the subcategory of tilting modules of $Rep(U_q\so_N)$,
see \cite{AP}. It only has finitely many simple objects up to isomorphism.
We only give some very basic facts about it here, see e.g. \cite{BK}
for more details.

We assume $q$ to be a primitive $2\ell$-th (for $N$ even) respectively
a primitive $4\ell$-th root of unity  (for $N$ odd), with $\ell \geq N$
in both cases.   The simple objects of $\F$ are labeled by 
the $\ell$-admissible  integral dominant
weights $\la$ of $\so_N$, which are defined by 
$$\la_1+\la_2 +N-2\leq \ell.$$
Theses categories also appear as level $\ell +2-N$ representations of
loop groups connected with $SO(N)$.
The image of the spinor representation $S$ is again simple and nonzero in
the quotient $\F$ and will be denoted by the same letter.
Tensoring a simpe object $V_\la$ by $S$ is given by the 
same rule \ref{tensorrule}, where we remove all representations
which are not labeled by an $\ell$-admissible weight.
With a little care, all the arguments we used above
will also go through for fusion categories related to $U_q\so_N$. 

\begin{theorem}\label{fusiontheorem} The statements of Theorem \ref{fft}
and \ref{sft} also hold in the context of fusion categories $\F$ as defined
in this section. In particular, $\End_\F(S^{\otimes n})$ is given by
a representation of $U'_q\o_n$ for $N$ even, and by a representation
of $U'_{-q^2}\so_n$ for $N$ odd.
\end{theorem}

$Proof.$ It is known that $S$ is a tilting module, and that $S^{\otimes n}$ can be written as
a direct sum of indecomposable tilting modules $T_\la$ with highest weight $\la$. 
The quotient of $S^{\otimes n}$ which lies in $\F$ is isomorphic to the summand consisting
of those $T_\la$ for which $\la$ is $\ell$-admissible. In this case $T_\la\cong V_\la$
is simple.
These basic facts allow us to proceed by induction on $n$ and on $N$ as in the proofs
for Theorem \ref{fft} and \ref{sft}. As there, one shows that $\End(S^{\otimes 2})$
is generated by $C$ (for $N$ odd) and by $C$ and $F$ (for $N$ even), using
the fact that the eigenvalues of $C$ for the different subrepresentations of $S^{\otimes 2}$
are  distinct as in the generic case (which would not necessarily be true if $q$ had a smaller degree).
Our conditions also ensure that the $q$-dimensions of all simple objects in $\F$
are nonzero. This is all needed to show that also the statement about
$S^{\otimes n}_{old}$ in the proof of Theorem \ref{fft} holds in this context.
We can now use the same argument as in the proof of Theorem \ref{fft} 
to show that $\End_\F(S^{\otimes n}_{new})$ will be isomorphic to a quotient
of $\End_{\F(U(N-2)}(S(N-2)_1^{\otimes n})$, where $\F(U(N-2))$ refers to
the fusion category for $U_q\so_{N-2}$ at the same root of unity.

\section{Related results and applications}

\subsection{Connection with results in \cite{Wspin} and Hasegawa duality}\label{Hasegawa}
This paper is closely related
to the paper \cite{Wspin}
which in turn was inspired by the paper \cite{Has} by Hasegawa. The key observation
there was (as far as this paper is concerned)  the fact that the commuting actions
of $O(N)$ and $SO(n)$ on $\C^N\otimes \C^n$ extend to commuting actions
of the corresponding spin groups on the Clifford algebra $Cl(Nn)\cong Cl(N)^{\otimes n}$
(isomorphism of vector spaces). This suggested commuting actions of these groups 
and the corresponding quantum groups on 
the $n$-th tensor power of the spin representation of $O(N)$. 
It was indeed shown in \cite{Wspin} for $N$ even that there exist actions 
of $\U=U_q\so_N\rtimes \Z/2$ and  and  $U'_q\so_n$
on $S^{\otimes n}$ which are each others commutants. As the Clifford algebra
is not simple for $N$ odd, no such simple statement could be shown in that case.
The best one could do was to prove a duality statement between the action of $U_q\so_N$ 
and the subalgebra of $U'_{q^2}\so_n$ generated by the elements $B_i^2$, $1\leq i<n$
(our $q^2$ 
here corresponds to $q$ in the parametrization in \cite{Wspin} for $N$ odd).
This subalgebra does not seem to allow a convenient algebraic description on its own.

In the current paper we do not apply Hasegawa's results at all. 
Using the operator $C$, we directly show that the centralizers can be
described in terms of 
the coideal algebra $U'_{-q}\so_n$. This does not
change much our previous description in terms of  $U'_{q}\so_n$ for
$N$ even, see Lemma \ref{so3rep}.
 But for $N$ odd, it gives the desired duality result
between actions of $U_q\so_N$ and  $U'_{-q^2}\so_n$, which now acts via
its non-classical representations on $S^{\otimes n}$.
Because of the latter, the result seems to be 
new even in the classical
case $q=1$. In particular, there does not seem to be any indication
of the algebra $U'_{-q}\so_n$ in the context of Hasegawa duality or
Howe duality (see next section).

\subsection{Connections to $q$-Howe duality}\label{qHowe} Similarly as for 
Clifford algebras in  Hasegawa duality, one obtains
commuting actions of groups, say $Gl(N)$ and $Gl(n)$ on
$\bigwedge (\C^N\otimes \C^n)\cong (\bigwedge \C^N)^{\otimes n}$.
For Lie type $A$, it has been shown  by Cautis, Kamnitzer and Morrison
\cite{CKM} that there similarly exist commuting actions of the quantum
groups $U_q\sl_N$ and $U_q\sl_n$ on $(\bigwedge \C^N)^{\otimes n}$.
The situation is more complicated for other Lie types, such as for
the commuting actions of
$\so_N$ and $\so_n$ as well as of $\sp_N$ and $\sp_n$ on the exterior
algebra $\bigwedge(\C^N\otimes \C^n)\cong (\bigwedge \C^N)^{\otimes n}$. 
 It was shown by Sartori and Tubbenhauer \cite{ST} that one obtains commuting 
actions of $U_q\so_N$ and $U'_q\so_n$ and of $U_q\sp_N$ and $U'_q\sp_n$,
where $U'_q\so_n$ is as in this paper, and $U'_q\sp_n\subset U_q\sl_n$
similarly is a coideal subalgebra. In the orthogonal case, at least for $N$
even, this result would also follow from 
 the results in \cite{Wspin}, as $S^{\otimes 2}\cong \bigwedge \C^N$
as $Pin(N)$-modules. The methods in \cite{ST} are completely different
from the ones in \cite{Wspin} and in this paper.

\subsection{Connections to the $q$-Clifford algebra in \cite{DF}}\label{DFconn}
A $q$-Clifford algebra has also been defined by Ding and Frenkel in \cite{DF}. Setting
$\hat\p_{k+1-i}=\Om_{i-1}\p_i$ and $\hat\pd_{k+1-i}=\Om_{i-1}\pd_i$,
one can show (see \cite{DF}, Proposition 5.3.1) that the elements
$\hat\p_i$ and $\hat\pd_i$ satisfy the relations of their $q$-Clifford algebra.
Similarly, if one defines $\hat\p_{i,-}$ and  $\hat\pd_{i,-}$
as before with $\Om_{i-1}$ replaced by $\Om_{i-1}^{-1}$, we obtain
the relations of their Clifford algebra with $q$ replaced by $q^{-1}$.
In particular the element $C$ would then have the somewhat more
appealing form
$$C=\sum_{i=1}^k \hat\p_i\otimes \hat\pd_{i,-}
+ \hat\pd_i\otimes \hat\p_{i,-}$$
in the even-dimensional case. The elements $\hat\p_i$, $\hat\p_{i,-}$ etc
already appeared in Section \ref{relationssec} as $c_{i,\pm}$ and $d_{i,\pm}$.
In particular, one can see the parameter $q$ appear explicitly in
the relations there, see e.g. \ref{cdrel}.

\subsection{Existence of representations of $U'_q\so_n$ and $U'_q\o_n$} The irreducible
representations of $U'_q\so_n$  for $q$ not a root 
of unity were constructed by Klimyk and his coauthors,
and in special cases also by different groups, see \cite{IK} and the references
there. Another approach was given by the author in \cite{Wspin}.
In all of these cases, the constructions of the representations were
somewhat  involved.  Our duality result provides yet another construction
of all irreducible representations of $U'_q\so_n$ and of many
irreducible representations of $U'_q\o_n$ if $q$ is not a root of unity,
see corollary \ref{onrep}. We also obtain some new
 nontrivial representations at roots of unity 
in the fusion cateogry setting. In our set-up, it suffices to find
the linear map $C$. After that the representations are obtained
in a comparatively  painless way using our duality result.

\subsection{New vertex models} Results in this paper were used by D. Gepner
in the recent preprint \cite{Gep} to construct new vertex models in connection with spin representations
and to study their algebraic properties.

\end{document}